\documentclass[11pt,a4paper]{article}
\usepackage{amssymb, setspace}
\usepackage{amsmath, eufrak,verbatim}
\usepackage{amsthm, graphicx, mathrsfs}
\usepackage[all]{xy}

\def\proof{\par\medskip\noindent {\sc Proof. }}

\def\intg           {\mathbb Z}
\def\rat            {\mathbb Q}
\def\rea            {\mathbb R}
\def\com            {\mathbb C}

\def\ra             {\rightarrow}

\def\hra            {\hookrightarrow}

\title{Moduli spaces of polarized irreducible symplectic manifolds are not necessarily connected}
\author{A. Apostolov}
\date{}
\begin{document}
\maketitle
\begin{abstract}
We show that the moduli space of polarized irreducible symplectic manifolds of $K3^{[n]}$-type, of fixed polarization type, is not
always connected. This can be derived as a consequence of Eyal Markman's characterization of polarized parallel-transport operators of $K3^{[n]}$-type.
\end{abstract}
\section{Introduction}
Irreducible (holomorphic) symplectic manifolds (also known as irreducible hyperk\"ahler manifolds) arise as elements in the Bogomolov decomposition of compact K\"ahler manifolds with trivial first Chern class (cf. [Bog]).
They can also be considered as higher-dimensional generalizations of the notion of a K3 surface.
A complex manifold $X$ is called an \emph{irreducible symplectic (IS) manifold} if $X$ is a simply-connected compact K\"ahler manifold such that
$H^0(X, \Omega_X^2)\cong \com \sigma$, where $\sigma$ is an everywhere non-degenerate holomorphic 2-form. The second cohomology
group $H^2(X,\intg)$ of an IS manifold carries an integral symmetric bilinear form $(\cdot,\cdot)_X$
of signature $(3,b_2(X)-3)$ called the \emph{Beauville form} (or \emph{Beauville-Bogomolov form}).
Another important invariant of IS manifolds is the \emph{Fujiki constant} $c_X\in \rat^+$.
It comes from the equality $\int_X \alpha^{2n}=c_X(\alpha,\alpha)^n_X, \forall \alpha\in H^2(X,\intg)$ (where $2n=\mathrm{dim}_{\com}X$).\\

  There are few known classes of examples of IS manifolds, namely deformations of
                                                                       \begin{itemize}
                                                                         \item Hilbert schemes of $n$ points on a K3 surface $S$, denoted by $S^{[n]}$;
these are known as IS manifolds \emph{of $K3^{[n]}$-type};
                                                                         \item generalized Kummer varieties, constructed as zero fibers of maps of the form $\Sigma \circ f$,
                                                                         where $f:\mathbb{T}^{[n]}\ra \mathbb{T}^{(n)}$ is the desingularization map from the Hilbert scheme of $n$ points
                                                                         on a complex torus $\mathbb{T}$ to the $n$-th symmetric power of the torus, and $\Sigma :\mathbb{T}^{(n)}\ra (\mathbb{T},0)$ is the sum map to the pointed torus $(\mathbb{T},0)$;
                                                                         \item O'Grady's examples in dimensions 6 and 10 (cf. [OG1] and [OG2]).
                                                                       \end{itemize}

For further details on the theory of IS manifolds, we refer to [Bea] and [Huy1].
There are several ways to organize IS manifolds into moduli spaces. In all cases one needs to add some extra structure to the data of an IS manifold.
For the theory of moduli spaces of IS manifolds, we follow the expositions in [Mar1, Ch. 1, 7] and [GHS2, Ch.3].
Fix a lattice $\Lambda$, which is isometric to $H^2(X,\intg)$ with its Beauville form $(\cdot,\cdot)_X$, for some IS manifold $X$.
A \textit{marked irreducible (holomorphic) symplectic manifold} $(X,\eta)$ consists of an IS manifold $X$ together with the choice of an isomorphism $\eta: H^2(X, \intg)\ra \Lambda$.
There is a (non-Hausdorff) coarse moduli space $\mathfrak{M}_{\Lambda}$, whose points represent equivalence classes $[(X,\eta)]$ of marked IS manifolds,
whose second cohomology group is isometric to $\Lambda$, where $(X,\eta)\sim (X',\eta')$ if there is an isomorphism $g:X\xrightarrow{\sim} X'$ such that $\eta'= g^*\circ \eta$ (cf. [Huy2, Ch. 3]).
Set
\begin{math}
\Omega_{\Lambda}=\{[w]\in \mathbb{P}(\Lambda\otimes_{\intg}\com)\ |\ (w,w)=0, (w, \overline{w})>0)\}.
\end{math}
There is a period map $P:\mathfrak{M}_{\Lambda}\ra \Omega_{\Lambda}$, which is a local homeomorphism ([Bea]) and is
surjective on each component of the moduli space ([Huy1, Thm. 8.1.]), sending $[(X,\eta)]$ to the period point $[\eta(\sigma)]\in \Omega_{\Lambda}$, where $\sigma$ is a generator of $H^0(X,\Omega_X^2)$.

We can also consider moduli spaces of (primitively) polarized IS manifolds. These have the advantage that their connected components admit the structure of
quasi-projective varieties. A \emph{polarization} on $X$ is the choice of an ample line bundle $L$ on $X$. This is equivalent to the choice of
a class $h=c_1(L)$ in $H^2(X,\intg)$, since the irregularity of $X$ is zero. Moduli spaces of polarized varieties with trivial canonical bundle were constructed by Viehweg ([Vie]),
as GIT quotients of Hilbert schemes. We can fix an IS manifold $X_0$ with $H^2(X_0, \intg)\cong \Lambda$ and the $O(\Lambda)$-orbit $\overline{h}$ of a primitive vector $h\in \Lambda$, of degree $(h,h)>0$.
This orbit $\overline{h}$ is called a \emph{polarization type}. Viehweg's construction yields a moduli space of polarized IS manifolds of type $(X_0,\overline{h})$ -- denote it by $\mathcal{V}_{X_0,\overline{h} }$. A point of $\mathcal{V}_{X_0, \overline{h}}$ represents an equivalence class of pairs $(X,L)$, where $X$ is a projective IS manifold, deformation equivalent to $X_0$, and $L\in Pic(X)$ is ample and such that $\eta(c_1(L))=h$ with respect to some marking $\eta$ of $X$. Given $h\in \Lambda$ with $(h,h)>0$,
denote by $\Omega_{h^{\perp}}$ the period domain $\Omega_{\Lambda}\cap h^{\perp}$ and let $O(\Lambda, h)$ be the stabilizer of $h$ in $O(\Lambda)$. This is a type IV symmetric domain and has two connected components - denote one of them by $\Omega_{h^{\perp}}^0$; let $O^+(\Lambda,h)$ be the subgroup of $O(\Lambda,h)$ fixing the connected components. Let $\mathcal{V}_{X_0, \overline{h}}^0$ be a component of $\mathcal{V}_{X_0, \overline{h}}$. There is a period map from
$\mathcal{V}_{X_0, \overline{h}}^0$ to $\Omega_{h^{\perp}}^0/O^+(\Lambda,h)$, which is a finite, dominant map of quasi-projective varieties ([GHS1, Thm. 1.5.]). In fact, even more is true: $\mathcal{V}_{X_0, \overline{h}}^0$ immerses as a dense open subset of a modular variety $\Omega_{h^{\perp}}^0/\Gamma_h$, where $\Gamma_h$ is a finite index
subgroup of $O^+(\Lambda, h)$ ([Mar1, Thm. 8.4.]).

In the case of K3 surfaces, fixing the degree of the polarization determines the polarization type as well -- the reason is that the K3 lattice is unimodular; in addition, the moduli space of
polarized K3 surfaces has one connected component for every degree of polarization. For IS manifolds of higher dimensions, it is natural to ask if this holds as well. First of all, it is easy to see that the degree does not determine the polarization type. For example, if $X_0$ is the Hilbert scheme of two points on a K3 surface, whenever the degree is congruent to -1 modulo 4, there are two polarization types for this degree -- \emph{split} and \emph{non-split} (cf. [GHS1, Ch. 3]). Then we can refine the question and ask if the moduli space of IS manifolds deformation equivalent to $X_0$, of fixed polarization type, is connected. The aim of this paper is to show that the answer is no in general - cf. Remark (1) at the end of Ch. 3 for explicit examples. More precisely, given positive integers $n, d$ and $t|\gcd(2n-2, 2d)$, define $\Sigma_n^{d,t}$ to be the set of isometry classes of pairs $(T,h)$, such that $T$ is an even positive definite lattice of rank two and discriminant $4d(n-1)/t^2$,  $h$ is a primitive element of square $(h,h)=2d$, and $h^{\perp}$ is generated by an element of square $2n-2$. In Ch. 2 we prove that there is a one-to-one correspondence between the set of connected components of the moduli space of polarized IS manifolds of $K3^{[n]}$-type, with polarization type of
degree $2d$ and divisibility $t$ and the set $\Sigma_n^{d,t}$. This turns out to be a consequence of Eyal Markman's characterization of polarized parallel-transport operators of $K3^{[n]}$-type, which yields an effective way of determining the number of connected components of the moduli space of polarized IS manifolds of $K3^{[n]}$-type. For the other known families of IS manifolds, such a characterization is still lacking but there are conjectural descriptions of their monodromy groups -- cf. the discussion in [Mar1, Ch. 10].\\

{\small \textbf{ACKNOWLEDGEMENTS}:  I would like to thank Eyal
Markman for explaining his work to me and Klaus Hulek and the anonymous referee for suggesting improvements to the exposition.}

\section{Monodromy Invariants}
In this section we give a short overview of the results from [Mar1] that we use - they describe a method for finding an invariant of the components of the
moduli space of polarized IS manifolds of $K3^{[n]}$-type. It
uses his characterization of parallel transport operators of $K3^{[n]}$-type (cf. [Mar1, Ch. 9.1.] and the references therein). The rough idea is to find a substructure in a certain lattice, whose isometry class remains unchanged under parallel transport.

First we introduce some notation. Let $\Lambda_n$ denote the lattice $E_8(-1)^{\oplus 2}\oplus U^{\oplus 3}\oplus \langle-2(n-1)\rangle$, and let
$\widetilde{\Lambda}$ denote the lattice $E_8(-1)^{\oplus 2}\oplus U^{\oplus 4}$, where $E_8(-1)$ is the rank 8 negative definite $E_8$ root lattice,
$U$ is the rank 2 hyperbolic lattice, and $\langle-2(n-1)\rangle$ is a rank one lattice generated by an element of length $-2(n-1)$. $\Lambda_n$
is an even lattice which is isometric to the Beauville lattice of an IS manifold of $K3^{[n]}$-type, whereas $\widetilde{\Lambda}$ is known as the \emph{Mukai lattice} (cf. [Muk]) --
it is denoted by $K(S)$ and is isometric to the lattice given by the
cohomology group of a $K3$-surface $H^*(S,\intg)=H^0(S,\intg)\oplus H^2(S,\intg)\oplus H^4(S,\intg)$, endowed with the \emph{Mukai pairing}
\begin{equation}
((r_1,h_1,s_1),(r_2,h_2,s_2))_M:=
(h_1,h_2)-r_1s_2-r_2s_1.
\end{equation}

From now on we only consider IS manifolds deformation equivalent to Hilbert schemes of points on a K3 surface. Put $X_0:=S^{[n]}$, where $S$ is a K3 surface.
There are finitely many connected components of $\mathfrak{M}_{\Lambda_n}$ parametrizing those marked pairs $(X,\eta)$, where $X$ is deformation equivalent to $X_0$ ([Mar1, L. 7.5.]).
Denote the set of these components by $\tau$ and let $\mathfrak{M}_{\Lambda_n}^{\tau}$ denote their union.
Let $\Omega_{\Lambda_n}$ and $\Omega_{h^{\perp}}$ denote the period domains associated to $\Lambda_n$ and $h^{\perp}\subset \Lambda_n$, introduced in the previous section.
There is an $O(\Lambda_n)$-equivariant \textit{refined period map}
$\widetilde{P}:\mathfrak{M}_{\Lambda_n}^{\tau}\ra \Omega_{\Lambda_n}\times \tau$, sending a marked pair $[(X,\eta)]$ to $([(X,\eta)], s)$, where
$\mathfrak{M}_{\Lambda_n}^{s}$ is the component containing $[(X,\eta)]$. $O(\Lambda_n)$ acts diagonally on $\Omega_{\Lambda_n}\times \tau$ by changing a period point
and a marking by an auto-isometry of $\Lambda_n$ (cf. [Mar1, Ch. 7.2.]).
Now $\Omega_{h^{\perp}}$ has two connected components -- the choice of $h$ and $s\in \tau$
determines the choice of a connected component of $\Omega_{h^{\perp}}$, which we denote by $\Omega_{h^{\perp}}^{s,+}$ (cf. [Mar1, Ch. 4]).
For $s\in \tau$, set $\mathfrak{M}_{h^{\perp}}^{s,+}:=\widetilde{P}^{-1}((\Omega_{h^{\perp}}^{s,+},s))$.
$\mathfrak{M}_{h^{\perp}}^{s,+}$ contains those marked pairs $[(X,\eta)]$ for which $\eta^{-1}(h)\in H^2(X,\intg)$ is of Hodge type $(1,1)$ and $(\eta^{-1}(h),\kappa)>0$ for a K\"ahler class $\kappa$, i.e. it is the first Chern class of a big line bundle (cf. [Huy1, Cor. 3.10.]).
Let $\mathfrak{M}_{h^{\perp}}^{s,a}\subset \mathfrak{M}_{h^{\perp}}^{s,+}$
be the subset, consisting of those $[(X,\eta)]$, for which $\eta^{-1}(h)$ is ample. $\mathfrak{M}_{h^{\perp}}^{s,a}$ is an open dense path-connected Hausdorff subset of $\mathfrak{M}_{h^{\perp}}^{s,+}$ ([Mar1, Cor. 7.3.]). This result uses the Global Torelli theorem of Verbitsky ([Ver, Thm. 1.16.]).
Let $\chi$ be an $O(\Lambda_n)$-orbit of pairs
$(h,s)$ with $(h,h)>0$ (this orbit is denoted by $\overline{h}$ in [Mar1]). Finally, form the disjoint union
\begin{equation*}
\mathfrak{M}_{\chi}^a:=\coprod_{(h,s)\in \chi}\mathfrak{M}_{h^{\perp}}^{s,a}.
\end{equation*}

$\mathfrak{M}_{\chi}^a$ is a coarse moduli space for marked polarized triples - a point of $\mathfrak{M}_{\chi}^a$ represents an isomorphism class $[(X,L,\eta)]$ of \emph{polarized marked IS manifolds of type} $(X_0,\chi)$.

Denote the moduli space of polarized IS manifolds of $K3^{[n]}$-type by $\mathcal{V}_{X_0}$ -- it is the union of the spaces $\mathcal{V}_{X_0, \overline{h}}$ over all polarization types $\overline{h}$, where $h\in \Lambda_n$ and $(h,h)>0$.
Let $[(X,H)]$ be a point in $\mathcal{V}_{X_0}$ and let $\mathcal{V}^0$ be the connected component containing $[(X,H)]$. The next proposition relates this component to a quotient of the coarse moduli space of marked polarized triples:\\

\textbf{Proposition 1.1. ([Mar1, Lemma 8.3.])}
\emph{There exists a natural isomorphism $\varphi: \mathcal{V}^0\ra \mathfrak{M}_{\chi}^a/O(\Lambda_n)$ in the category of analytic spaces.}\\

To any family $\mathcal{X}\ra T$ of IS manifolds such that there is a point $t_0\in T$ with $\mathcal{X}_{t_0}\cong X$, one can associate a natural monodromy representation
$\pi_1(T,t_0)\ra GL(H^*(X,\intg))$. A \emph{monodromy operator} is an isomorphism of the cohomology $H^*(X, \intg)$ in the image of some monodromy representation.
The subgroup of $GL(H^*(X,\intg))$ generated by all
monodromy operators is denoted by $\mathrm{Mon}(X)$. $\mathrm{Mon}^2(X)$ denotes the image of $\mathrm{Mon}(X)$ in $GL(H^2(X,\intg))$. In fact, since monodromy operators preserve the
Beauville-Bogomolov pairing, $\mathrm{Mon}^2(X)$ is a subgroup of $O(H^2(X,\intg))$.

Let $O(\Lambda_n,\widetilde {\Lambda})$ (resp. $O(H^2(X,\intg),\widetilde {\Lambda})$) denote the set of primitive isometric embeddings of $\Lambda_n$ (resp. $H^2(X,\intg)$) into $\widetilde{\Lambda}$. Now any IS manifold $X$ determines an orientation class $or_X\in H^2(\widetilde{\mathcal{C}}_X,\intg)$,  i.e. a
generator of $H^2(\widetilde{\mathcal{C}}_X,\intg)\cong \intg$, where $\widetilde{\mathcal{C}}_X$ is the cone
$\{h\in H^2(X,\rea)|\ (h,h)>0\}$ (cf. [Mar1, Ch. 4]). The cone $\widetilde{\mathcal{C}}_{\Lambda_n}$ in $\Lambda_n \otimes \rea$ is defined analogously and the set of
its orientations is denoted by $\mathrm{Orient}(\Lambda_n)$ -- it is the set of two generators of $H^2(\widetilde{\mathcal{C}}_{\Lambda_n}, \intg)$. $O^+(H^2(X, \intg))$ (resp. $O^+(\Lambda_n)$)
denotes the subgroup of $O(H^2(X, \intg))$ (resp. $O(\Lambda_n)$) whose elements act trivially on $H^2(\widetilde{\mathcal{C}}_X,\intg)$ (resp. $H^2(\widetilde{\mathcal{C}}_{\Lambda_n}, \intg)$).

By studying the representation of the
monodromy group on the cohomology of $X$, E. Markman came up with the following idea -- let $X$ be of $K3^{[n]}$-type and assume first that $n\geq 4$. Now consider $Q^4(X,\intg)$, which is the quotient of $H^4(X,\intg)$ by the image of the cup product homomorphism
$\cup: H^2(X, \intg)\otimes H^2(X, \intg)\ra H^4(X,\intg)$.
Now, $Q^4(X,\intg)$ admits a monodromy invariant bilinear pairing which makes it isometric to the Mukai lattice. Moreover, the $\mathrm{Mon}(X)$-module $\mathrm{Hom}[H^2(X,\intg),Q^4(X,\intg)]$contains a unique rank 1 saturated $\mathrm{Mon}(X)$
-submodule, which is a sub-Hodge structure of type $(1,1)$ ([Mar1, Thm. 9.3.]). A generator of this module induces an $O(\Lambda_n)$-orbit of primitive isometric embeddings of $H^2(X,\intg)$
into the Mukai lattice, such that the image of $H^2(X,\intg)$ under such an embedding is orthogonal to the image of the projection of $c_2(X)\in H^4(X,\intg)$ in $Q^4(X,\intg)$. As for the case
$n = 2, 3$ - there is only a single $O(\widetilde{\Lambda})$-orbit of primitive isometric embeddings of $H^2(X, \intg)$ in the Mukai lattice $\widetilde{\Lambda}$ anyway.
This yields the following statement:\\

\textbf{Theorem 1.2. ([Mar1, Cor. 9.5.])}
\emph{Let $X$ be an IS manifold of $K3^{[n]}$-type, $n\geq 2$. $X$ comes with a natural choice of an $O(\widetilde{\Lambda})$-orbit $[\iota_X]$ of primitive isometric
embeddings of $H^2(X, \intg)$ in the Mukai lattice $\widetilde{\Lambda}$. The subgroup $Mon^2(X)$ of $O^+(H^2(X, \intg))$ is equal to the stabilizer of $[\iota_X]$ as an element of the orbit space
$O(H^2(X, \intg),\widetilde{\Lambda})/O(\widetilde{\Lambda})$.}\\

\textbf{Proposition 1.3. (cf. [Mar1, Cor. 9.10.])}
\emph{The set $\tau$ of connected components of the moduli space of marked IS manifolds of $K3^{[n]}$-type is in bijective correspondence to the orbit set
$[O(\Lambda_n,\widetilde {\Lambda})/O(\widetilde{\Lambda})]\times \mathrm{Orient}(\Lambda_n)$, where
$O(\widetilde{\Lambda})$ acts by post-composition on $O(\Lambda_n,\widetilde {\Lambda})$.}\\

The bijection is given by mapping a component $s$ to the pair $([\iota_X\circ \eta^{-1}],\eta_*(or_X))$, where $[(X,\eta)]$ is a point of $\mathfrak{M}_{\Lambda_n}^s$.

Next we introduce parallel-transport operators in the polarized setting:

\textbf{Definition 1.4. ([Mar1, Def. 1.1.(4)])}
\emph{Let $X_1, X_2$ be IS manifolds and let $H_i\in \mathrm{Pic}(X_i)$ be ample. Set $h_i:=c_1(H_i)$. An isomorphism $f:H^2(X_1,\intg)\ra H^2(X_2,\intg)$ is said to be a
\textbf{polarized parallel-transport operator} from $(X,H_1)$ to $(X,H_2)$, if there exists a smooth, proper family of
IS manifolds $\pi:\mathcal{X}\ra T$ onto an analytic space, points $t_1,t_2\in T$, isomorphisms $\psi_i:X_i\ra \mathcal{X}_{t_i}$, $i=1,2$, a continuous path $\gamma: [0,1]\ra T$
with $\gamma(0)=t_1, \gamma(1)=t_2$, and a flat section $h$ of $R^2\pi_*\intg$, such that parallel transport in the local system $R^2\pi_*\intg$ along $\gamma$ induces the
homomorphism $(\psi_{2}^{-1})^*\circ f\circ \psi_1^*: H^2(\mathcal{X}_{t_1}, \intg)\ra H^2(\mathcal{X}_{t_2}, \intg)$, $h_{t_i}=(\psi_{i}^{-1})^*(h_i)$, $i=1,2$, and $h_t$
is an ample class in $H^{1,1}(\mathcal{X}_t, \intg)$, $\forall t\in T$.}\\

 We call an isometry $g:H^2(X_1, \intg)\ra H^2(X_2, \intg)$ \emph{orientation-preserving}, if $g^*or_{X_2}=or_{X_1}$.
We can now state the following characterization of polarized parallel-transport operators:\\

\textbf{Theorem 1.5. (cf. [Mar1, Cor. 7.4., Thm. 9.8.])}
\emph{Let $(X_1,H_1)$ and $(X_2,H_2)$ be polarized IS manifolds of $K3^{[n]}$-type. Set $h_i:=c_1(H_i), i=1,2$. An isometry $g:H^2(X_1, \intg)\ra H^2(X_2, \intg)$ is a polarized parallel-transport operator from $(X_1,H_1)$ to $(X_2,H_2)$ if and only if
$g$ is orientation-preserving,
\begin{math}
[\iota_{X_1}]=[\iota_{X_2}]\circ g,
\end{math}
and $g(h_1)=h_2$.}\\

The above statement can be used to obtain a lattice-theoretic characterization of polarized parallel-transport operators in the following manner -- choose a primitive isometric embedding $\iota_X: H^2(X, \intg)\hra \widetilde{\Lambda}$ in the $O(\widetilde{\Lambda})$-orbit given by Thm 1.2.
For a primitive class $h\in H^2(X, \intg)$, of degree $(h,h)=2d>0$, let $T(X,h)$ denote the saturation in $\widetilde{\Lambda}$, of the sublattice spanned by
$\iota_X(h)$ and $\mathrm{Im}(\iota_X)^{\perp}$.
$T(X,h)$ is a rank 2 positive definite lattice.
Denote by $[(T(X,h), \iota_X(h))]$ the isometry class of the pair $(T(X,h), \iota_X(h))$, i.e.
$(T(X,h), \iota_X(h))\sim (T', h')$ iff there exists an isometry $\gamma:T(X,h)\ra T'$ such that $\gamma(\iota_X(h))=h'$.
Let $I(X)$ denote the set of primitive cohomology classes of positive degree in $H^2(X, \intg)$.
Let $\Sigma_n$ be the set of isometry classes of pairs $(T,h)$, consisting of an even rank 2 positive definite lattice $T$ and a primitive element $h\in T$
such that $h^{\perp}\cong \langle 2n-2\rangle$.
Let $f_{X}: I(X)\ra \Sigma_n$ be the function sending $h$ to $[(T(X,h), \iota_X(h))]$.
Note that $[(T(X,h), \iota_X(h))]$ does
not depend on the choice of representative of the orbit $[\iota_X]$.
Also $f_X(h)=f_{X'}(h')$, for any isomorphism $X\xrightarrow{\sim} X'$ mapping $h'$ to $h$ in cohomology. $f_{X}$ is called
a \textit{faithful monodromy invariant function} in [Mar2] because it separates orbits for the action of the monodromy group of $X$ on $I(X)$ (cf. [Mar2, Ch. 5.3.]).\\

\textbf{Proposition 1.6. ([Mar3, Lemma 0.4.])}
\emph{Let $(X_1, H_1)$ and $(X_2, H_2)$ be two polarized pairs of
IS manifolds of $K3^{[n]}$-type. Set $c_1(H_i)=h_i$. Then $f_{X_1}(h_1)=f_{X_2}(h_2)$ if and only if there exists a polarized parallel-transport operator from $(X_1, H_1)$ to $(X_2, H_2)$.}
\proof One direction is clear -- suppose there exists a polarized parallel-transport operator $g: H^2(X_1, \intg)\ra H^2(X_2, \intg)$
from $(X_1, H_1)$ to $(X_2, H_2)$. In particular, $g(h_1)=h_2$. Since, by Thm. 1.5.,
$[\iota_{X_1}]=[\iota_{X_2}]\circ g$, there exists an isometry $\gamma\in O(\widetilde{\Lambda})$ such that $\gamma\circ \iota_{X_1}=\iota_{X_2}\circ g$.
The map $\gamma$ induces an isometry between the pairs $(T(X_1,h), \iota_{X_1}(h_1))$ and $(T(X_2,h),\iota_{X_2}(h_2))$, i.e. $f_{X_1}(h_1)=f_{X_2}(h_2)$.

Now assume that $f_{X_1}(h_1)=f_{X_2}(h_2)$. This means that the two pairs $(T(X_1,h), \iota_{X_1}(h_1))$ and $(T(X_2,h),\iota_{X_2}(h_2))$ are isometric.
The idea is to construct the required parallel transport operator $g$ from an isometry of $\widetilde{\Lambda}$.
Now, $T(X_1,h)$ and $T(X_2,h)$ are primitively-embedded sublattices of signature $(2, 0)$, of the same isometry class, in the unimodular lattice
$\widetilde{\Lambda}$, of signature $(4, 20)$. Therefore, [Nik, Thm. 1.1.2.b)] implies that there exists a $\gamma\in O(\widetilde{\Lambda})$, such that $\gamma(T(X_1,h))=T(X_2,h)$,
and $\gamma(\iota_{X_1}(h_1))=\iota_{X_2}(h_2)$. Set $g=\iota_{X_2}^{-1}\circ \gamma \circ \iota_{X_1}$:
\[\xymatrix{H^2(X_1, \intg)\ar^{\ \ \ \  \ \ \iota_{X_1}}@^{(->}[r]\ar^{g}[d] & \widetilde{\Lambda}\ar^{\gamma}[d] & T(X_1,h)\ar^{\gamma|_{T(X_1,h)}}[d]\ar@{*{\supset\ \ \ \ \ \ }{}}[l]\\
 H^2(X_2, \intg)\ar^{\ \ \ \  \ \ \iota_{X_2}}@^{(->}[r] & \widetilde{\Lambda} & T(X_2,h)\ar@{*{\supset\ \ \ \ \ \ }{}}[l]}\]

Indeed, the map $g$ maps $H^2(X_1, \intg)$ to $H^2(X_2, \intg)$ because 
\begin{equation*}
\gamma(\iota_{X_1}(H^2(X_1, \intg))^{\perp})=\iota_{X_2}(H^2(X_2, \intg))^{\perp}.
\end{equation*}
Then 
\begin{equation*}
[\iota_{X_2}]\circ g=[\iota_{X_2}\circ (\iota_{X_2}^{-1}\circ \gamma \circ \iota_{X_1})]=[\gamma \circ \iota_{X_1}]=[\iota_{X_1}],
\end{equation*} 
 since $\gamma\in O(\widetilde{\Lambda})$. Furthermore, 
 \begin{equation*}
 g(h_1)=\iota_{X_2}^{-1}\circ \gamma \circ \iota_{X_1}(h_1)=\iota_{X_2}^{-1}\circ \iota_{X_2} (h_2)=h_2.
 \end{equation*}
Now assume that $g$ is orientation-reversing. Choose $\alpha\in  H^2(X_2, \intg)$ satisfying $(\alpha, \alpha)=2$, $(\alpha, h_2)=0$. Define the isometry
$\rho_{\alpha}(\lambda):=-\lambda + (\alpha, \lambda)\alpha $, $\lambda\in H^2(X_2, \intg)$. Set $\tilde{g}:=-\rho_{\alpha}\circ g$. Since
$\rho_{\alpha}$ is an element of $\mathrm{Mon}^2(X_2)$ ([Mar1, Thm. 9.1.]) and $\rho_{\alpha}(h_2)=-h_2$, $\tilde{g}$ is an orientation-preserving isometry between $H^2(X_1, \intg)$ and $H^2(X_2, \intg)$, satisfying
$[\iota_{X_1}]=[\iota_{X_2}]\circ \tilde{g}$ and $\tilde{g}(h_1)=h_2$.
Thm. 1.5. implies that $\tilde{g}$ is a polarized parallel-transport operator from $(X_1, H_1)$ to $(X_2, H_2)$. \hfill $\blacksquare$\\

\section{Enumerating the components}

Let $\mu_n$ denote the set of connected components of $\mathcal{V}_{X_0}$. We are now ready to prove the following\\

\textbf{Theorem 2.1.}
\emph{There is an injective map $f:\mu_n \ra \Sigma_n$, given by mapping a connected component $\{\mathcal{V}^0\}$ of $\mathcal{V}_{X_0}$ to $f_{X}(h)$ for some $[(X,H)]\in \mathcal{V}^0$.}
\proof
First of all, the map $f$ is well-defined.
Pick a point $[(X_1,H_1)]\in \mathcal{V}^0$. Choose a marking $\eta_1$ of $X_1$ and let $\mathfrak{M}_{h^{\perp}}^{a,s_1}$ be the component of the moduli space of
marked polarized triples, containing the point $[(X_1, H_1, \eta_1)]$. Choose another point $[(X_2,H_2)]\in \mathcal{V}^0$ and a marking $\eta_2'$ of $X_2$. Then
$[(X_2, H_2,\eta_2')]\in \mathfrak{M}_{\chi}^{a}$, where $\chi$ is the
$O(\Lambda_n)$-orbit of $(h,s_1)$. Since
$O(\Lambda_n)$ acts transitively on the set of components of $\mathfrak{M}_{\chi}^a$, there is a marking $\eta_2$ such that
$[(X_2, H_2, \eta_2)]\in \mathfrak{M}_{h^{\perp}}^{a,s_1}$. Choose a path $\gamma:[0,1]\ra \mathfrak{M}_{h^{\perp}}^{a,s_1}$ such that
$\gamma(0)=[(X_1, H_1, \eta_1)]$, $\gamma(1)=[(X_2, H_2, \eta_2)]$. For each $r\in [0,1]$, $\gamma(r)$ has a Kuranishi neighborhood $U_r\subset \mathfrak{M}_{h^{\perp}}^{a,s_1}$ and a
semi-universal family of deformations $\pi_r:\mathcal{X}_r \ra U_r$. Upon shrinking $U_r$, we may assume that it is simply-connected. As in the proof of [Mar1, Cor. 7.4.], we can use these to construct a polarized parallel transport operator from $(X,H_1)$ to $(X,H_2)$, namely $\gamma([0,1])$ admits a finite subcovering $\{U_{r_i}\}_{i\in \{1,...,m\}}$
and we can choose a partition $s_0=0, s_1, ..., s_{m-1}, s_m=1$ of $[0,1]$ such that $\gamma(s_i)\in U_{r_i}\cap U_{r_{i+1}}, i=1,...,m-1$ and $\gamma([s_{i-1},s_i])\subset U_{r_i}$.
Now set $V_1:=U_{r_1}$, $V_i:=V_{i-1}\bigcup_{\gamma(s_{i-1})}U_{r_{i}}$ for $i=2,...,m$, where $V_{i-1}\bigcup_{\gamma(s_{i-1})}U_{r_{i}}$ denotes the pushout of the inclusions of the point $\gamma(s_{i-1})$ in $V_{i-1}$ and $U_{r_{i}}$. We can identify the fibers $\pi_{r_i}^{-1}(\gamma(s_i))$ and $\pi_{r_{i+1}}^{-1}(\gamma(s_i))$ in order to obtain a family $\pi: \mathcal{X} \ra V_m$ of IS manifolds, as well as a path $\tilde{\gamma}:[0,1]\ra V_m$ which is the concatenation of the paths $\gamma([s_{i-1},s_i]), i=1,...,m-1$,
and which yields a polarized parallel-transport operator from $(X_1, H_1)$ to $(X_2, H_2)$. Hence, by Prop. 1.6., $f_{X_1}(h_1)=f_{X_2}(h_2)$, and the map $f$ is well-defined.

Now let $\mathcal{V}^1$ and $\mathcal{V}^2$  be two components of $\mathcal{V}_{X_0}$ such that $f(\{\mathcal{V}^1\})=f(\{\mathcal{V}^2\})$. Choose points $[(X_i,H_i)]\in \mathcal{V}^i$, $i=1,2$.
Since $f_{X_1}(h_1)=f_{X_2}(h_2)$, there exists a polarized parallel-transport operator from $(X_1, H_1)$ to $(X_2, H_2)$. Indeed, by Prop. 1.5. and 1.6., there is an analytic family of IS manifolds
$\pi: \mathcal{X} \ra T$ and a section $\tilde{h}$ of $R^2\pi_*\intg$ such that $(\mathcal{X}_{t_0},\tilde{h}_{t_0})\cong (X_1,h_1)$ and $(\mathcal{X}_{t_1},\tilde{h}_{t_1})\cong (X_2,h_2)$ for some $t_1,t_2\in T$.
Consider the path $\gamma: [0,1]\ra T$ with $\gamma(0)=t_1, \gamma(1)=t_2$ (cf. Def. 1.4.). The long exact sequence associated to the exponential sequence yields:
\begin{equation*}
...\ra R^1\pi_*\mathcal{O}_{\mathcal{X}}\ra R^1\pi_*\mathcal{O}_{\mathcal{X}}^*\ra R^2\pi_*\intg \ra R^2\pi_*\mathcal{O}_{\mathcal{X}} \ra...
\end{equation*}
Since $\tilde{h}_t$ has type (1,1) at every point $t \in T$, the image of $\tilde{h}$ in $H^0(\mathcal{X}, R^2\pi_*\mathcal{O}_{\mathcal{X}})$ vanishes.
The coherent sheaf $R^1\pi_*\mathcal{O}_{\mathcal{X}}$ also vanishes,
since the irregularity of the fibers is zero. Hence $h$ lifts locally to sections of the relative analytic Picard sheaf $R^1\pi_*\mathcal{O}_{\mathcal{X}}^*$.
Then we may find a finite open cover $\{W_i\}_{i=1}^m$ of $\gamma([0,1])$ by simply connected sets $W_i$, markings
$\tilde{\eta}^i$, and bundles $\mathcal{H}_i\in \mathrm{Pic}(\mathcal{X}|_{W_i}/W_i)$ which make $(\mathcal{X}|_{W_i}, \mathcal{H}_i, \tilde{\eta}^i)$ into families of marked polarized triples, such that the classifying maps
$\alpha_i$ of the families map $W_i$ to the same component $\mathfrak{M}_{h^{\perp}}^{a,0}$ of the moduli space of marked polarized triples $\mathfrak{M}_{\chi}^{a}$,
for all $i=1,..., m$. Let $q:\mathfrak{M}_{\chi}^{a}\ra \mathfrak{M}_{\chi}^a/O(\Lambda_n)$ be
the quotient map. Prop. 1.1. gives analytic isomorphisms $\varphi_i: \mathcal{V}^i\ra \mathfrak{M}_{\chi}^a/O(\Lambda_n), i=1,2$.
Consider the points $[(X_2',H_2')]:=\varphi_1^{-1}(q(\alpha_m(t_1)))\in \mathcal{V}^1$ and $[(X_2,H_2)]\in \mathcal{V}^2$.
Since they map to the same orbit of marked polarized triples in $\mathfrak{M}_{\chi}^a/O(\Lambda_n)$, there is an analytic isomorphism between $(X_2',H_2')$ and $(X_2,H_2)$.
By GAGA (cf. [JPS]), it is induced by an algebraic isomorphism between $(X_2',H_2')$ and $(X_2,H_2)$.
 Hence $[(X_2',H_2')]=[(X_2,H_2)]$ as points in $\mathcal{V}^{\tilde{\tau}}$, which implies that $\mathcal{V}^1=\mathcal{V}^2$. \hfill $\blacksquare$

Fix a K3 surface $S$ and an isometry $(H^*(S,\intg),(\cdot,\cdot)_M)\cong \widetilde{\Lambda}$ (cf. (1)). Let $X$ be of $K3^{[n]}$-type and let $h_d\in H^2(X, \intg)$ be a primitive element of degree $2d$ and choose $\iota$ in the orbit $[\iota_X]$, so that $\mathrm{Im}(\iota)^{\perp}$ is generated by the vector
$v:=(1,0,1-n)\in \widetilde{\Lambda}$. Let $\langle 2d\rangle$ denote a rank 1 lattice generated by an element of length $2d$.
For  $r, s\in \intg, (r,s)$ denotes the gcd of $r$ and $s$. Set $t:=\mathrm{div}(h_d)$ to be the \emph{divisibility} of $h_d$, i.e. the positive generator of the ideal $(h_d, H^2(X,\intg))$ in $\intg$. In particular, $\intg/t \intg$ is a subgroup of the discriminant group of $H^2(X,\intg)$ and the latter has order $2n-2$. This means that $t$ divides $2n-2$.
The next proposition relates $t$ to the index of
$\langle \iota(h_d)\rangle \oplus \langle v\rangle$ in the saturation $T(X,h_d)$ with respect to $\iota_X$. I thank the referee for suggesting a simpler proof of this proposition.\\

\textbf{Proposition 2.2.}

\emph{The integer $t$ is equal to the index of $\langle \iota(h_d)\rangle \oplus \langle v\rangle$ in $T(X,h_d)$.}
\proof
$\iota(h_d)$ can be written as $(c,tm\xi, c(n-1))$ in $\widetilde{\Lambda}$, where $c, m\in \intg, \xi \in H^2(S, \intg)$ is primitive and
$(m, \frac{2n-2}{t})=1$. Moreover, $(c,m)=1$ by the primitivity of $\iota(h_d)$ in $v^{\perp}$. Now consider the class $u:=\frac{\iota(h_d)-cv}{t}=(0,m\xi, \frac{c(2n-2)}{t})$.
It is integral because of $t\ |\ 2n-2$, and it is primitive in $\widetilde{\Lambda}$ because of $(m,\frac{c(2n-2)}{t})=1$. Then it is easy to check that the lattice $\langle u, v\rangle$ is
saturated in $\widetilde{\Lambda}$. First of all, pick a class $w:=av+bu=(a, bm\xi, (1-n)(a-b\frac{2c}{t}))$ which is primitive in $\langle u, v\rangle$; this means that $(a,b)=1$. Suppose that
$w=f\tilde{w}$, where $\tilde{w}$ is primitive in $\widetilde{\Lambda}$. Then $f$ divides $a$, $bm$ and $b\frac{c(2n-2)}{t}$. Furthermore, since $(a,b)=1$, $f$ divides $(m, \frac{c(2n-2)}{t})=1$,
i.e. $f=1$ and $w=\tilde{w}$. This implies that $\langle u, v\rangle$ is saturated, i.e. $T(X,h_d)=\langle u, v\rangle$. Finally, since the index of $\langle \iota(h_d)\rangle \oplus \langle v\rangle$ in $\langle u, v\rangle$ is equal to $t$ by a direct discriminant computation, we have shown the claim. \hfill $\blacksquare$\\

Now given positive integers $d$ and $t|(2d,2n-2)$, \emph{let $\mu^{d,t}_n\subset \mu_n$} be the subset of connected components of $\mathcal{V}_{X_0}$
parametrizing polarized IS
manifolds of $K3^{[n]}$-type, with polarization type of
degree $2d$ and divisibility $t$. By [GHS1, Cor. 3.7.], this data is sufficient to determine the polarization type whenever $(\frac{2n-2}{t},\frac{2d}{t})$ and $t$ are coprime.

\textbf{Proposition 2.3.}
\emph{The image of the restriction $f|_{\mu^{d,t}_n}$ is $\Sigma_n^{d,t}$.}

\proof

Let $(T,h)$ be a representative of some isometry class in $\Sigma_n^{d,t}$. Let $\delta$ generate $h^{\perp}$ in $T$.
Then, by [Nik, Thm. 1.1.2.] we may choose a primitive isometric embedding
$j:T\hra \widetilde{\Lambda}$. Let $\tilde{\iota}: \Lambda_n\hra \widetilde{\Lambda}$ be a primitive isometric embedding. Then the sublattice $\tilde{\iota}(\Lambda_n)^{\perp}$
is generated by a primitive element of degree $2n-2$ and we may find an embedding $\iota$ in the $O(\widetilde{\Lambda})$-orbit
$[\tilde{\iota}]$ that $\iota(\Lambda_n)^{\perp}=\langle j(\delta) \rangle$,
since $O(\widetilde{\Lambda})$ acts transitively on primitive elements. In particular, $j(h)\in \iota(\Lambda_n)$. The embedding $\iota$ determines a pair of
components $t_1,t_2 \in \tau$ of the moduli space $\mathfrak{M}^{\tau}_{\Lambda_n}$ of marked IS manifolds of $K3^{[n]}$-type, via the bijection $\tau\cong [O(\Lambda_n,\widetilde{\Lambda})/O(\widetilde{\Lambda})]\times \mathrm{Orient}(\Lambda_n)$ (cf. Prop. 1.3.).
Set $h_1:=\iota^{-1}(j(h))\in \Lambda_n$. Since  $\mathfrak{M}_{h_1^{\perp}}^{a,t_1}$ is nonempty,
we may choose a point in it, represented by a marked IS manifold $(X,\eta)$. Then, by the definition of $\mathfrak{M}_{h_1^{\perp}}^{a,t_1}$, $h_2:=\eta^{-1}(h_1)$
is an ample class in $H^2(X,\intg)$. Let $\chi$ be the $O(\Lambda_n)$-orbit of $(h_1,t_1)$ and let $\mathcal{V}^0\cong \mathfrak{M}_{\chi}^a/O(\Lambda_n)$ be the component of the moduli space of polarized pairs, containing $[(X,H_2)]$. Then, by construction,
$f(\{\mathcal{V}^0\})=[(T,h)]$, i.e. $f|_{\mu^{d,t}_n}$ is surjective onto $\Sigma_n^{d,t}$.\hfill $\blacksquare$\\

We obtain as a corollary:

\textbf{Corollary 2.4.}
\emph{The number of connected components of the moduli space of polarized IS manifolds of $K3^{[n]}$-type, with polarization type of
degree $2d$ and divisibility $t$, is given by $|\Sigma_n^{d,t}|$}.

\proof

Injectivity of $f|_{\mu^{d,t}_n}$ is given by Thm. 2.1; surjectivity onto $\Sigma_n^{d,t}$,by Prop. 2.3. Hence the number of connected components, i.e. the
number of elements in the image of $f|_{\mu^{d,t}_n}$, is given by$|\Sigma_n^{d,t}|$. \hfill $\blacksquare$\\

\section{Computations}

In this section, we compute the cardinality of the set $\Sigma_n^{d,t}$ and give examples.
For an integer $r$, $\varphi (r)$ denotes the Euler $\varphi$-function, and $\rho(r)$ denotes the
number of prime divisors of $r$.\\

\textbf{Proposition 3.1.}
\begin{it}

Let $t$ be a divisor of $(2d, 2n-2)$, and set $D:=4d(n-1)/t^2$, $g:=(2d, 2n-2)/t$, $\tilde{n}:=(2n-2)/(2d, 2n-2), \tilde{d}:=2d/(2d, 2n-2)$, $w:=(g, t)$, $g_1:=g/w$,
$t_1:=t/w$. Put $w=w_+(t_1)w_-(t_1)$ where $w_+(t_1)$ is the product of all powers of primes dividing $(w, t_1)$. Then

     \begin{itemize}
       \item $|\Sigma_n^{d,t}|=w_+(t_1)\varphi(w_-(t_1))2^{\rho(t_1)-1}$ if $t>2$ and one of the following sets of conditions hold: $g_1$ is even, $(\tilde{d}, t_1)=(\tilde{n}, t_1)=1$ and
the residue class $-\tilde{d}/\tilde{n}$ is a quadratic residue modulo $t_1$; \textbf{OR} $g_1, t_1$, and $\tilde{d}$ are odd, $(\tilde{d}, t_1)=(\tilde{n}, 2t_1)=1$ and
$-\tilde{d}/\tilde{n}$ is a quadratic residue modulo $2t_1$; \textbf{OR} $g_1, t_1$, and $w$ are odd, $\tilde{d}$ is even, $(\tilde{d}, t_1)=(\tilde{n}, 2t_1)=1$ and
$-\tilde{d}/(4\tilde{n})$ is a quadratic residue modulo $t_1$;
       \item $|\Sigma_n^{d,t}|=w_+(t_1)\varphi(w_-(t_1))2^{\rho(t_1/2)-1}$ if $t>2$, $g_1$ is odd, $t_1$ is even, $(\tilde{d}, t_1)=(\tilde{n}, 2t_1)=1$ and
$-\tilde{d}/\tilde{n}$ is a quadratic residue modulo $2t_1$;
       \item $|\Sigma_n^{d,t}|=1$ if $t\leq 2$ and one of the following sets of conditions hold: $g_1$ is even, $(\tilde{d}, t_1)=(\tilde{n}, t_1)=1$ and
the residue class $-\tilde{d}/\tilde{n}$ is a quadratic residue modulo $t_1$; \textbf{OR} $g_1, t_1$, and $\tilde{d}$ are odd, $(\tilde{d}, t_1)=(\tilde{n}, 2t_1)=1$ and
$-\tilde{d}/\tilde{n}$ is a quadratic residue modulo $2t_1$; \textbf{OR} $g_1, t_1$, and $w$ are odd, $\tilde{d}$ is even, $(\tilde{d}, t_1)=(\tilde{n}, 2t_1)=1$ and
$-\tilde{d}/(4\tilde{n})$ is a quadratic residue modulo $t_1$; \textbf{OR} $g_1$ is odd, $t_1$ is even, $(\tilde{d}, t_1)=(\tilde{n}, 2t_1)=1$ and
$-\tilde{d}/\tilde{n}$ is a quadratic residue modulo $2t_1$;
       \item $|\Sigma_n^{d,t}|=0$, else.
     \end{itemize}
\end{it}
\proof
Denote by $D(T)$ the discriminant group of a lattice $T$. By definition, $|\Sigma_n^{d,t}|$ is equal to the number of isometry classes of primitive inclusions of $\langle 2d\rangle$
into lattices $T$ of discriminant $D$, such that the orthogonal complement of $\langle 2d\rangle$ in a given $T$ is isometric to $\langle 2n-2\rangle$. By [Nik, Thm. 1.5.1], the latter
are classified by certain equivalence classes of monomorphisms
$\gamma:H\ra D(\langle 2n-2\rangle)\cong \intg/ (2n-2)\intg$, where $H$ is the unique subgroup of $D(\langle 2d\rangle)\cong \intg/2d\intg$
of order $t$ such that $\Gamma_{\gamma}$ (the graph of $\gamma$)
is an isotropic subgroup of $D(\langle 2d\rangle)\oplus D(\langle 2n-2\rangle)$ with respect to the discriminant form. The equivalence is given by: $\gamma\sim \gamma'$ if there exist $\psi\in O(\langle 2n-2\rangle)\cong \intg/2\intg, \varphi \in O(\langle 2d\rangle)\cong \intg/2\intg$
such that $\gamma\circ \overline{\varphi}=\overline{\psi}\circ \gamma'$ ($\overline{\varphi}$ and $\overline{\psi}$ are the induced isometries on the discriminant groups).
Now the set of subgroups of the form $\Gamma_{\gamma}$ can be identified with the set of generators in $D(\langle 2d\rangle)\oplus D(\langle 2n-2\rangle)$ of the form
$(\tilde{d}g, \tilde{n}g\varepsilon)$, where $\varepsilon$ is a unit modulo the index $t$. The isotropy condition on $\Gamma_{\gamma}$ reads as follows:
\begin{equation*}
\frac{(\tilde{d}g)^2}{2d}+\frac{(\tilde{n}g\varepsilon)^2}{2n-2}
=\frac{\tilde{d}g+\tilde{n}g\varepsilon^2}{t}= \frac{\tilde{d}g_1+\tilde{n}g_1\varepsilon^2}{t_1}\equiv 0\ (\mathrm{mod\ } 2\intg).
\end{equation*}

In other words, we seek the number of classes $\varepsilon$ modulo $t$, coprime to $t$, which satisfy the congruence $\tilde{d}+\tilde{n}\varepsilon^2\equiv 0\ (\mathrm{mod}\ 2t_1)$.
The conditions for existence as well as the number of such solutions to a congruence of this form were already determined in the proof of
Prop. 3.6 of [GHS1] in the context of computing orbits of vectors under the stable orthogonal group.
Now, $O(\langle 2d\rangle)$ and $O(\langle 2n-2\rangle)$ act on $\Gamma_{\gamma}$ by 'reflecting', i.e. by flipping signs in the first, resp. second coordinate of the elements of the
graph. In addition, $\Gamma_{\gamma}$ is central-symmetric (i.e. $(a,b)\in \Gamma_{\gamma}\Leftrightarrow (-a,-b)\in \Gamma_{\gamma}$). For $2\geq t$, the graphs are fixed by the
action, for $t>2$ there are no fixed graphs, hence, by the above considerations
there are two graphs in each equivalence class, so we need to divide the numbers from [GHS1, Prop. 3.6] by two and we obtain the result. \hfill $\blacksquare$\\

In the following cases the polarization type is determined by the values of $t$ and $d$, and the corresponding moduli spaces are indeed connected:

\textbf{Proposition 3.2.}
\emph{The moduli space of polarized IS manifolds of $K3^{[n]}$-type is connected in each of the following cases:}.\\
 \begin{itemize}
    \item $t=1$, any degree $2d>0$ and dimension $2n>2$ such that $(n-1,d)=1$ ('split' polarization);
    \item $t=2$, any degree $2d>0$ and dimension $2n>2$ such that $(n-1,d)=1$ and $d+n-1\equiv 0 (\mathrm{mod\ } 4)$('non-split' polarization, cf. [GHS1, Ch.3]);
    \item $t=p^{\alpha}$ for a prime $p>2$, any degree $2d$ and dimension $2n>2$ such that $p^{2\alpha}|(2d,2n-2)$.
     \end{itemize}

\textbf{Remarks.}

(1) \textbf{Example:} In particular, Prop. 3.2. implies that the moduli space of polarized IS manifolds of $K3^{[2]}$-type, with fixed polarization type, is connected.
However, the following example shows that fixing the polarization type need not imply connectedness of the moduli space.
Let $d=pq$ and $n-1=mpq$, where $p$ and $q$ are different primes, and $-m$ is a quadratic residue modulo $pq$. Set $t=pq$.
In this case, the polarization type is determined by $t$ and $d$, since $(2m,2)=2$ and $t=pq$ are coprime. But $|\Sigma_n^{d,t}|=2$, i.e. there are
two connected components with this polarization type.

(2) Whenever the moduli space of IS manifolds of fixed polarization type has more than one component, we can pick two points $[(X_1,H_1)]$ and $[(X_2,H_2)]$ from different components;
suppose that they map to the same period point. This yields a Hodge isometry $H^2(X_1, \intg)\ra H^2(X_2, \intg)$, which is not a polarized parallel-transport operator from $[(X_1,H_1)]$ to $[(X_2,H_2)]$. In particular, there need not exist a birational map from $X_1$ to $X_2$. In an e-mail to the author, E. Markman raised the interesting question of constructing an algebraic correspondence in $X_1\times X_2$, which induces this isometry, for every such pair of points in the moduli space.

(3) Note that the formulas in Prop.2.1. imply that
the number of connected components $|\Sigma_n^{d,t}|$ can get arbitrarily large, as we
vary the dimension $2n-2$ and the degree $2d$.

\end{document}